\documentclass[10pt,journal,twoside,web]{ieeecolor}

\usepackage{generic,cite}
\usepackage{amsmath,amssymb,amsfonts,mathtools}
\usepackage{multirow,array}
\usepackage[spaces,hyphens]{xurl} % allow arbitrary line breaks in a formatted URL string
\usepackage[colorlinks,allcolors=blue]{hyperref}

\usepackage[ruled,linesnumbered]{algorithm2e}
\SetKw{Continue}{continue}
\SetKw{Break}{break}
\SetKw{Initialize}{Initialize}
\SetKwComment{Comment}{/* }{ */}
\SetKwComment{SingleComment}{// }{}
\SetKwComment{TriComment}{$\triangleright$\ }{}

\numberwithin{equation}{section}
\DeclareMathOperator*{\argmin}{arg\,min}

\newtheorem{theorem}{Theorem}[section]
\newtheorem{lemma}[theorem]{Lemma}

\newtheorem{definition}{Definition}[section]
\newtheorem{remark}{Remark}[section]

\renewcommand{\theequation}{\arabic{section}.\arabic{equation}}

\def\BibTeX{{\rm B\kern-.05em{\sc i\kern-.025em b}\kern-.08em
    T\kern-.1667em\lower.7ex\hbox{E}\kern-.125emX}}
\begin{document}

\title{Projection-based Prediction-Correction Method for Distributed Consensus Optimization}
\author{Han Long
\thanks{Han Long is with the Southern University of Science and Technology, 1088 Xueyuan Avenue, Shenzhen 518055, P.R. China (e-mail: 12232848@sustech.edu.cn).}
}

\maketitle

\begin{abstract}
Within the realm of industrial technology, optimization methods play a pivotal role and are extensively applied across various sectors, including transportation engineering, robotics, and machine learning. With the surge in data volumes, there is an increasing demand for solving large-scale problems, which in turn has spurred the development of distributed optimization methods. These methods rely on the collaborative efforts of numerous dispersed devices to achieve the collective goals of the system. This study focuses on the exploration of distributed consensus optimization problems with convex set constraints within networks. The paper introduces a novel Adaptive Projection Prediction-Correction Method (PPCM), inspired by the proximal point algorithm and incorporating the theory of variational inequalities. As a contraction algorithm with notable convergence performance, PPCM is particularly suited for decentralized network environments. Moreover, the selection of parameters for this method is both straightforward and intuitive, avoiding the complexities of intricate parameter tuning. Comprehensive theoretical analysis and empirical testing have validated the effectiveness of PPCM. When applied to problems such as distributed linear least squares, logistic regression, and support vector machines, PPCM demonstrates superior performance, achieving computation speeds over ten times faster than built-in Python functions while maintaining high precision. In conclusion, this research provides a valuable distributed consensus optimization technique, both theoretically and practically.
\end{abstract}

\section{Introduction}
\label{sec:intro}
\par In recent years, mathematical optimization algorithms have emerged as indispensable components of industrial technology development and applications. Optimization plays a pivotal role across diverse domains, including transportation engineering\cite{su2009traffic,mohebifard2018distributed,teodorovic2008swarm}, robot control\cite{chen2018cooperative,tron2016distributed,li2013contact}, smart grids\cite{zhang2012convergence,braun2016distributed,yu2018economic}, machine learning\cite{joshi2022distributed} and more.

\par As data continues to accumulate, the scale of problems has expanded. The traditional single-centralized computing model no longer meets the demands of engineering. Simultaneously, the rapid progress in electronic systems, chip industries, and communication technologies has facilitated the establishment of network systems to tackle intricate problems. Distributed computing, with its collaborative nature within network systems, effectively handles the ever-increasing computational requirements\cite{qiu2016survey,sun2019survey}. Consequently, distributed optimization has gained prominence and significance in decision-making and data processing across various fields. Successful applications encompass energy economic dispatch\cite{yang2016distributed}, distributed control\cite{nedic2018distributed}, robot systems\cite{lynch2017modern}, sensor networks\cite{rabbat2004distributed}, smart buildings\cite{minoli2017iot}, intelligent manufacturing\cite{davis2015smart}, and beyond.

\par Distributed optimization involves the cooperative collaboration of dispersed intelligent devices or microprocessors to achieve system-level objectives. In the context of engineering systems, distributed optimization empowers subsystems to make localized decisions while interacting with each other to pursue optimal system performance. In computer science, distributed optimization is frequently employed to distribute computationally intensive training tasks across multiple microprocessors and coordinate their efforts towards a coherent training goal. Although specific circumstances may vary, the fundamental concept of distributed optimization entails decomposing a comprehensive mathematical optimization problem into smaller-scale subproblems and authorizing multiple computing agents to solve these subproblems in a coordinated manner, thereby approximating the optimal or near-optimal solution to the original mathematical optimization problem.

\par Distributed consensus optimization represents a crucial research domain within distributed systems and networks, aiming to address how multiple computational entities (e.g., servers, sensors, or robots) can achieve a consensus decision. From a mathematical perspective, distributed consensus optimization fundamentally pertains to a class of optimization problems where the objective function exhibits finite sum and structured characteristics, typically manifested as
\begin{equation}
    \min\limits_{x \in \mathcal{X}} f(x) := \sum\limits_{i=1}^N f_i(x).
    \label{finite-sum}
\end{equation}
Numerous practical applications feature problems that exhibit a resemblance or equivalence in form of \eqref{finite-sum}. For instance, the least square solution of linear equations arises when estimating the electromechanical oscillation modes of extensive power system networks using synchrophasors\cite{nabavi2015distributed}. In the domain of sensor networks, various application scenarios can be identified, such as robust estimation\cite{zielinski1983pj}, energy-based source localization\cite{chen2002source,sheng2003energy}, and distributed expectation maximization\cite{nowak2003distributed}. Additionally, comparable problems manifest in the realm of large-scale machine learning\cite{abeynanda2021study}. There are many other practical applications, which will not be exhaustively enumerated here.

\par Currently, a multitude of researchers have devised a wide array of algorithms to address distributed consensus optimization problems. Among the notable contributions, Nedi\'{c} \& Ozdaglar\cite{nedic2009distributed,nedic201010} introduced a method predicated on weighted averaging for subgradient approaches in scenarios devoid of constraints, subsequently extending their framework to incorporate convex set limitations\cite{nedic2010constrained}. Another innovative approach was presented by Zanella et al.\cite{varagnolo2015newton}, who unveiled a consensus-based method leveraging the Newton-Raphson algorithm. The inception of employing the push-sum consensus model for devising distributed optimization strategies was attributed to Tsianos et al.\cite{tsianos2012consensus}, a methodology that has garnered extensive exploration and refinement in various studies thereafter\cite{tsianos2011distributed,tsianos2012push,tsianos2013role}. Furthermore, Pu et al.\cite{pu2020push} developed a novel push-pull gradient technique. This method distinguishes itself by a dual mechanism where gradient information is disseminated to adjacent nodes (push), whilst decision variable information is assimilated from them (pull), thereby coining the term ``push–pull gradient methods.” Additionally, the consensus ADMM\cite{chang2014multi} strategy emerged, notable for its employment of an inexact step during each ADMM update. This strategic choice permits the execution of computationally economical operations at every iteration, enhancing the algorithm's efficiency and applicability in distributed settings.

\par Despite the plethora of algorithms available for addressing distributed consensus issues, each exhibiting a unique blend of strengths and weaknesses, this paper introduces a novel decentralized algorithm that primarily relies on gradients and projections. Characterized by its simplicity and clarity in parameter selection, this algorithm is readily adaptable to a diverse array of complex and dynamic distributed network environments. It is underpinned by a rigorous yet accessible theoretical framework for convergence, demonstrating rapid convergence rates. Empirical evaluations have showcased its exceptional performance, positioning it as a potent new tool for tackling distributed consensus optimization challenges.

\section{Preliminaries}
\label{sec:pre}

\subsection{Norm}
\label{subsec:norm}
\begin{definition}[H-norm of a vector]
Let $H$ be a symmetric positive-definite matrix. Then, the H-norm of a vector $x$ is defined as
\begin{equation*}
    \|x\|_H = \sqrt{x^T H x}.
\end{equation*}
\end{definition}

\begin{definition}[H-norm of a matrix]
Let $H$ be a symmetric positive-definite matrix. The H-norm of a matrix $A$ is defined using the induced norm, that is
\begin{equation*}
    \|A\|_H = \sup\limits_{x \neq 0} \frac{\|Ax\|_H}{\|x\|_H}.
\end{equation*}
\end{definition}

\begin{lemma}
\label{lemma_Hnorm}
Let $H$ be a symmetric positive-definite matrix. The H-norm of a matrix $A$ is equal to its 2-norm. In other words,
\begin{equation*}
    \|A\|_H = \sqrt{\rho(A^T A)} = \|A\|
\end{equation*}
\end{lemma}
\begin{proof}
\begin{align*}
    \|A\|_H &= \sup\limits_{x \neq 0} \frac{\|Ax\|_H}{\|x\|_H} = \sup\limits_{x \neq 0} \sqrt{\frac{x^T A^T H A x}{x^T H x}} \\
    &= \sup\limits_{y \neq 0} \sqrt{\frac{y^T H^{-\frac{1}{2}} A^T H A H^{-\frac{1}{2}} y}{y^T y}} \qquad \bigg(\mathrm{let}\ y = H^{\frac{1}{2}}x\bigg) \\
    &= \sup\limits_{y \neq 0} \frac{\|H^{\frac{1}{2}} A H^{-\frac{1}{2}} y\|}{\|y\|} = \| H^{\frac{1}{2}} A H^{-\frac{1}{2}} \| \\
    &= \sqrt{\rho\bigg(H^{\frac{1}{2}} A^T A H^{-\frac{1}{2}}\bigg)} = \sqrt{\rho(A^T A)} \\
    &= \| A \|.
\end{align*}
\end{proof}

\par In the definitions and lemma above, we assume that matrix (and vector) multiplication is dimensionally compatible. Furthermore, throughout this paper, $\| \cdot \|$ denotes the Euclidean norm (or 2-norm).

\subsection{Graph theory}
\label{sec_graph}
\par Denote an graph $\mathcal{G} = (\mathcal{V}, \mathcal{E})$ with $N$ vertexes, where $\mathcal{V} = \{1, \cdots, N\}$ is the node set and $\mathcal{E} \subseteq \mathcal{V} \times \mathcal{V}$ is the edge set. Define the neighbour of $i$ by $\mathcal{N}_i = \{ j \in \mathcal{V} | (i, j) \in \mathcal{E} \}$ and $d_i = |\mathcal{N}_i|$ is the degree of $i$. The graph is undirected if whenever $(i,j) \in \mathcal{E}$ implies that $(j,i) \in \mathcal{E}$. An undirected graph $\mathcal{G}$ is said to be connected if for any two nodes $i, j$ there exists a path from $i$ to $j$, i.e. one can find $v_0 = i, v_1, \cdots, v_m = j$ such that $(v_{k-1}, v_k) \in \mathcal{E}$ for $1 \leq k \leq m$. 
\par The adjacency matrix of $\mathcal{G}$ is defined by $A = (a_{ij})_{N \times N}$ with $a_{ij} > 0$ if $(i,j) \in \mathcal{E}$ and $a_{ij} = 0$ otherwise. Then $A = A^T$ when $\mathcal{G}$ is undirected. Assume that the graph does not have any multiedges and self loops(i.e. $a_{ii} = 0,\ \forall i \in \mathcal{V}$).
\par The Laplacian matrix of the undigraph $\mathcal{G}$ is defined by $L = (l_{ij})_{N \times N}$ with $l_{ii} = \sum\limits_{j \neq i} a_{ij}$ and $l_{ij} = -a_{ij}$ for all $i \neq j$. Note that $L$ is symmetric and has zero row sums. Then $0$ is an eigenvalue of $L$ with the corresponding eigenvector $e := (1, \cdots, 1)^T$, the $N \times 1$ column vector of ones. Also, $L$ is diagonally dominant with nonnegative diagonal entries, thus it is positive semidefinite. When $\mathcal{G}$ is connected, $0$ is a simple eigenvalue of $L$. Let $\lambda_i(L)$ be the $i$-th smallest eigenvalue of $L$ which gives that $0 = \lambda_1(L) \leq \lambda_2(L) \leq \cdots \leq \lambda_N(L)$. Besides, $\lambda_2(L)$ is called the algebraic connectivity, which is positive if and only if $\mathcal{G}$ is connected\cite{merris1994laplacian}. Additionally, by utilizing the Gershgorin circle theorem, we can obtain $\rho(L) \leq 2 \| A \|_{\infty}$.

\subsection{Properties of projection}
\label{sec_proj}

\par Let $\Omega \subseteq \mathbb{R}^n$ be a closed, convex set. The projection onto $\Omega$ under the Euclidean norm $P_{\Omega}(\cdot)$ is defined by
\begin{equation*}
    P_{\Omega}(v) = \argmin\limits_{u \in \Omega} \|u-v\|,\quad v \in \mathbb{R}^n.
\end{equation*}

\begin{lemma}
Let $\Omega \subseteq \mathbb{R}^n$ be closed and convex. Then
\begin{equation}
    (v - P_{\Omega}(v))^T(u - P_{\Omega}(v)) \leq 0,\quad \forall v \in \mathbb{R}^n,\ \forall u \in \Omega.
    \label{eq_proj_2}
\end{equation}
\label{lemma_proj}
\end{lemma}
\begin{proof}
By the definition of projection, the following inequality holds:
\begin{equation}
    \|v - P_{\Omega}(v)\| \leq \|v - w\|,\quad \forall v \in \mathbb{R}^n,\ \forall w \in \Omega.
    \label{eq_proj_3}
\end{equation}
Notice that for any $v \in \mathbb{R}^n$, $P_{\Omega}(v) \in \Omega$. Since $\Omega$ is closed and convex, then for any $u \in \Omega$ and $\theta \in (0,1)$, we have
\begin{equation*}
    w := \theta u + (1 - \theta) P_{\Omega}(v) = P_{\Omega}(v) + \theta(u - P_{\Omega}(v)) \in \Omega.
\end{equation*}
Take this $w$ in \eqref{eq_proj_3}, we obtain
\begin{equation*}
    \|v - P_{\Omega}(v)\|^2 \leq \|v - P_{\Omega}(v) - \theta(u - P_{\Omega}(v))\|^2.
\end{equation*}
Expanding the above inequality yields
\begin{equation*}
    (v - P_{\Omega}(v))^T(u - P_{\Omega}(v)) \leq \frac{\theta}{2}\|u - P_{\Omega}(v)\|^2.
\end{equation*}
Let $\theta \rightarrow 0_+$ and the lemma is proved.
\end{proof}

\subsection{Optimal condition in variation inequality form}
\label{sec:vi}
\par Consider the general convex optimization problem with linear constraint
\begin{equation}
    \begin{aligned}
        \min\limits_{x}\quad & f(x) \\
        \text{s.t.}\quad & Cx = d, \\
        & x \in \mathcal{X},
    \end{aligned}
    \label{prob_coplc}
\end{equation}
where $f:\mathbb{R}^n \rightarrow \mathbb{R}$ is convex and differentiable, $\mathcal{X} \subseteq \mathbb{R}^n$ is closed and convex, $C \in \mathbb{R}^{m \times n}$, $d \in \mathbb{R}^m$. By introducing the Lagrangian multiplier $\lambda \in \Lambda = \mathbb{R}^m$ the linear constraint $Cx-d=0$, the Lagrange function of \eqref{prob_coplc} is
\begin{equation}
    L(x, \lambda) = f(x) - \lambda^T(Cx - d),\quad x \in \mathcal{X},\ \lambda \in \Lambda.
    \label{eq_lagrange}
\end{equation}
Assuming that all the regularity conditions are met. Therefore, there exists a saddle point $(x^*, \lambda^*) \in \mathcal{X} \times \Lambda$ which satisfies
\begin{equation}
    L_{\lambda \in \Lambda}(x^*, \lambda) \leq L(x^*, \lambda^*) \leq L_{x \in \mathcal{X}}(x, \lambda^*).
    \label{eq_saddle_point}
\end{equation}
Here $x^*$ is a solution of \eqref{prob_coplc}. The equivalent form of \eqref{eq_saddle_point} is
{\footnotesize
\begin{equation}
    \left\{
    \begin{aligned}
        &x^* \in \mathcal{X}, &f(x) - f(x^*) + (x - x^*)^T(-C^T \lambda^*) &\geq 0, &\forall x \in \mathcal{X}, \\
        &\lambda^* \in \Lambda, &\hfill (\lambda - \lambda^*)^T(Cx^* - d) &\geq 0, &\forall \lambda \in \Lambda.
    \end{aligned}
    \right.
    \label{eq_vi_saddle_point}
\end{equation}
}
Utilizing the properties of differentiable convex functions, we have
\begin{equation}
    f(x) - f(y) \geq (x - y)^T \nabla f(y),\quad \forall x \in \mathcal{X},\ \forall y \in \mathcal{X}
    \label{eq_convex_1}
\end{equation}
and
\begin{equation}
    (x - y)^T [\nabla f(x) - \nabla f(y)] \geq 0,\quad \forall x \in \mathcal{X},\ \forall y \in \mathcal{X}.
    \label{eq_convex_2}
\end{equation}
If we denote
\begin{equation*}
    u = \left(
    \begin{array}{c}
        x \\
        \lambda
    \end{array}
    \right),\quad
    F(u) = \left(
    \begin{array}{c}
        g(x) - C^T\lambda \\
        Cx - d
    \end{array}
    \right),\quad
    \Omega = \mathcal{X} \times \Lambda,
\end{equation*}
where $g = \nabla f$ denotes the gradient of $f$. Apply \eqref{eq_convex_1} to \eqref{eq_vi_saddle_point} after substituting $x$ with $(1-\theta)x^* + \theta x$ and let $\theta \rightarrow 0^+$, then \eqref{eq_saddle_point} can be written in the variational inequality form:
{\small
\begin{flalign}
    \text{VI}(\Omega, F) &&
    u^* \in \Omega^*,\ (u - u^*)^T F(u^*) \geq 0,\ \forall u \in \Omega. &&
    \label{eq_vi}
\end{flalign}
}
Here $u^*$ represents a solution of VI$(\Omega, F)$, and conventionally, we denote the solution set as $\Omega^*$. Notice that
{\footnotesize
\begin{equation*}
    F(u) = \left(
    \begin{array}{c}
        g(x) - C^T\lambda \\
        Cx - d
    \end{array}
    \right) = 
    \left(
    \begin{array}{cc}
        0 & -C^T \\
        C & 0
    \end{array}
    \right)
    \left(
    \begin{array}{c}
        x \\
        \lambda
    \end{array}
    \right)
    +
    \left(
    \begin{array}{c}
        g(x) \\
        -d
    \end{array}
    \right),
\end{equation*}
}
and with \eqref{eq_convex_2}, we obtain
\begin{flalign*}
    &\forall u \in \Omega,\ \forall \tilde{u} \in \Omega, \\
    &\quad (u - \tilde{u})^T [F(u) - F(\tilde{u})] = (x - \tilde{x})^T [g(x) - g(\tilde{x})] \geq 0.
\end{flalign*}
Then we can deduce that $F$ is monotone, thereby implying that
\begin{equation}
    u^* \in \Omega^*,\quad (u - u^*)^T[F(u) - F(u^*)] \geq 0,\quad \forall u \in \Omega.
    \label{eq_monotone_2}
\end{equation}

\begin{remark}
If $f$ is a convex function but not differentiable, we only need to introduce the subgradient by taking $ g \in \partial f $. Utilizing the properties of the subgradient, the same conclusions can be drawn without any modifications.
\end{remark}

\section{The Distributed Optimization Model}
\par We focus on the distributed consensus optimization problem of the following form
\begin{flalign}
    \text{(P)} &&
    \begin{aligned}
    \min\limits_{x}\quad & \sum\limits_{i=1}^N f_i(x) \\
    \text{s.t.}\quad & x \in \mathcal{X} = \cap_{i=1}^N \mathcal{X}_i
    \end{aligned}
    &&
    \label{P}
\end{flalign}
in a network of $N$ agents, where $\mathcal{X}_i \subseteq \mathbb{R}^n$ is a closed, convex set which is known to agent $i$. Each $f_i : \mathbb{R}^n \rightarrow \mathbb{R}$ is the private objective function kept only by agent $i$. We assume that all the $f_i$'s are convex and differentiable, thus (P) is a convex optimization problem.
\par The network can be mathematically described by a connected undirected graph $\mathcal{G} = (\mathcal{V}, \mathcal{E})$ with $N$ vertexes. Edge $(i, j) \in \mathcal{E}$ indicates that agent $i$ and $j$ are able to communicate with each other. One direct reformulation of (P) that accounts for the networked communication is
\begin{flalign}
    \text{(P1)} &&
    \begin{aligned}
        \min\limits_{(x_i)_{i=1}^N}\quad & \sum\limits_{i=1}^N f_i(x_i) \\
        \text{s.t.}\quad & x_i = x_j,\ \forall (i, j) \in \mathcal{E}; \\
        & x_i \in \mathcal{X}_i,\ 1 \leq i \leq N.
    \end{aligned}
    &&
    \label{P1}
\end{flalign}
\par To obtain a more compact form, let
\begin{equation*}
    x = \left(
    \begin{array}{c}
        x_1 \\
        \vdots \\
        x_N
    \end{array}
    \right) \in \mathbb{R}^{n \times N},\quad
    f(x) = \sum\limits_{i=1}^N f_i(x_i),\quad
    C = L \otimes I_n,
\end{equation*}
here $L$ is the Laplacian matrix of $\mathcal{G}$, $\otimes$ denotes the Kronecker product and $I_n$ is the identity matrix of order $n$. Then we can reformulate (P1) by
\begin{flalign}
    \text{(P2)} &&
    \begin{aligned}
        \min\limits_{x}\quad & f(x) \\
        \text{s.t.}\quad & Cx = 0, \\
        & x \in \mathring{\mathcal{X}} := \mathcal{X}_1 \times \cdots \times \mathcal{X}_N,
    \end{aligned}
    &&
    \label{P2}
\end{flalign}
which is a convex optimization problem with linear constraint.

\begin{remark}
In this work, we represent communication networks with undirected graphs, given that bidirectional communication is common in real-world networks, while unidirectional communication is largely an artificial and uncommon setup. We translate the network communication need into a linear constraint $Cx = 0$, equivalent to the consensus constraint where all $x_i$ are equal. Thus, our proposed algorithm can adapt to and handle both directed and dynamically changing networks, provided this equivalence is met.
\end{remark}

% \subsection{Related work}
% \par Currently, a multitude of researchers have devised a wide array of algorithms to address problems analogous to the form delineated in (P). Nedi\'{c} \& Ozdaglar\cite{nedic2009distributed,nedic201010} proposed the weighted-averaging-based subgradient method for unconstrained cases, and later, Nedi\'{c} et al.\cite{nedic2010constrained} further refined the method to accommodate convex set constraints. Zanella et al.\cite{varagnolo2015newton} introduced a Newton-Raphson consensus-based method. Tsianos et al.\cite{tsianos2012consensus} were the pioneers in utilizing the push-sum consensus model to develop distributed optimization methods, which have since been extensively explored in subsequent works\cite{tsianos2011distributed,tsianos2012push,tsianos2013role}. Pu et al.\cite{pu2020push} devised a push-pull gradient method. Additionally, several algorithms are grounded in the alternating direction method of multipliers (ADMM), such as classical ADMM\cite{wei2012distributed}, consensus ADMM\cite{chang2014multi}, flexible ADMM\cite{hong2016convergence}, edge-based flexible ADMM\cite{wei20131}, asynchronous ADMM\cite{zhang2014asynchronous,chang2016asynchronousp1,chang2016asynchronousp2}, and more.

\section{The Projection-based Prediction-Correction Method}
\par Recall that we have transformed (P) into the form of (P2). According to section \ref{sec:vi}, it suffices to solve the variational inequality VI$(\Omega, F)$. In this section, we introduce our novel method and provide a rigorous convergence analysis, inspired by the approach presented in the seminal work by He\cite{he1997class}.

\subsection{Motivation}
\par Presently, a variety of techniques exist for tackling problems in the form of \eqref{prob_coplc}. Nevertheless, effectively managing convex constraints remains a formidable task. Specifically, within the ADMM framework, when confronted with box constraints, resolving the subproblems can be notably intricate, at times even rivaling or surpassing the intricacy of the initial quandary. Given real-world scenarios where acquiring gradients is typically viable and the process of projecting onto convex sets is often computationally efficient, there is a predilection for favoring gradient-driven methodologies complemented by projection strategies for effectively addressing such challenges. A comprehensive and elegant proof of the following lemma can be found in \cite{he1999inexact}.
\begin{lemma}
$u^*$ is a solution of VI$(\Omega, F)$ if and only if $u^*$ is a solution of the projection equation
\begin{equation}
    u = P_{\Omega}[u - \beta F(u)],
    \label{eq_vipe}
\end{equation}
where $\beta > 0$.
\label{thm_vipe}
\end{lemma}
\par Now, let us proceed to solve VI$(\Omega, F)$. One viable method is the proximal point algorithm (PPA), which generates the iteration sequence $\{u^k\}$ satisfying the following condition:
{\footnotesize
\begin{flalign*}
    &u^{k+1} \in \Omega, \notag \\
    &\quad (u - u^{k+1})^T[F(u^{k+1}) + \frac{1}{\beta_k}(u^{k+1} - u^k)] \geq 0,\quad \forall u \in \Omega,\ \forall k \in \mathbb{N}.
\end{flalign*}
}
With Theorem \ref{thm_vipe}, we can concisely express the iterations of PPA by
\begin{equation}
    u^{k+1} = P_{\Omega}[u^k - \beta_k F(u^{k+1})].
    \label{it_ppa}
\end{equation}
Thought PPA has nice convergence property, \eqref{it_ppa} clearly represents an implicit iterative scheme, posing significant challenges in solving each iteration.
\par Another feasible approach is the extra-gradient method proposed by Korpelevich\cite{korplevich1976ekstragradientnyi}. It can be formulated in a prediction-correction manner as follows:
\begin{equation}
    \left\{
    \begin{aligned}
        &\tilde{u}^k = P_{\Omega}[u^k - \beta F(u^k)], \\
        &u^{k+1} = P_{\Omega}[u^k - \beta F(\tilde{u}^k)].
    \end{aligned}
    \right.
    \label{method_eg}
\end{equation}
The extra-gradient method converges if $\beta$ is chosen appropriately. Also, He\cite{he2002improvements} has enhanced the extra-gradient method by
\begin{equation}
    \left\{
    \begin{aligned}
        &\tilde{u}^k = P_{\Omega}[u^k - \beta_k F(u^k)], \\
        &u^{k+1} = P_{\Omega}[u^k - \alpha_k\beta_k F(\tilde{u}^k)].
    \end{aligned}
    \right.
    \label{method_he_2}
\end{equation}
Here, $\beta_k$ only needs to satisfy the local Lipschitz condition, and $\alpha_k$ can be computed in an adaptive manner. Furthermore, He\cite{he2004comparison} provides the following method:
\begin{equation}
    \left\{
    \begin{aligned}
        &\tilde{u}^k = P_{\Omega}[u^k - \beta_k F(u^k)], \\
        &d(u^k, \tilde{u}^k) := (u^k - \tilde{u}^k) - \beta_k[F(u^k) - F(\tilde{u}^k)], \\
        &u^{k+1} = u^k - \alpha_k d(u^k, \tilde{u}^k).
    \end{aligned}
    \right.
    \label{method_he_1}
\end{equation}
The projection operation is omitted in the correction step. Essentially, method \eqref{method_he_2} and method \eqref{method_he_1} offer two distinct descent directions. Coincidentally, these two directions can be iterated using the same step size $\alpha_k$, earning them the term ``twin directions". The detailed proof can be found in \cite{he1992globally}.
\par However, the above methods overlooks the separable structure of $F$ in the variational inequality. By fully exploiting this characteristic, we contemplate the development of a semi-explicit scheme, more precisely
\begin{equation}
    \left\{
    \begin{aligned}
        x^{k+1} &= P_{\mathcal{X}}\big[x^k - \beta_k\big(g(x^k) - C^T \lambda^{k}\big)\big], \\
        \lambda^{k+1} &= P_{\Lambda}[\lambda^k - \beta_k(Cx^{k+1} - d)].
    \end{aligned}
    \right.
    \label{it_semiex}
\end{equation}
Although the aforementioned scheme might not guarantee convergence, it presents a promising avenue. Leveraging the iteration yielded by \eqref{it_semiex} as a prediction, we then implement the correction step to guide the iterative sequence towards the solution point, ensuring eventual convergence.
\par Revisiting (P2), the intricacies of decentralized distributed systems pose a formidable challenge when it comes to updating iterations through the utilization of a uniform step-size, say $\beta_k$. In an ideal setting, each node autonomously determines a fitting step-size parameter tailored to its specific requirements. Hence, our algorithm is crafted by drawing inspiration from He's research\cite{he2018uniform}, and its design is outlined as follows in response to this inherent complexity.

\subsection{Scheme of the proposed method}
\par Let
\begin{equation*}
    R_k = \left(
    \begin{array}{ccc}
        r_1 I_n & & \\
        & \ddots & \\
        & & r_N I_n
    \end{array}
    \right), \
    S_k = \left(
    \begin{array}{ccc}
        s_1 I_n & & \\
        & \ddots & \\
        & & s_N I_n
    \end{array}
    \right)
\end{equation*}
be positive definite and block diagonal. Define
\begin{equation}
    H_k = \left(
    \begin{array}{cc}
        R_k & \\
        & S_k^{-1}
    \end{array}
    \right).
    \notag
\end{equation}
Let $\eta \in (0,1)$ and $\mu$ be constant such that $\|C^T\| \leq \mu$. Also, let $g$ denote the gradient of $f$ defined in \eqref{P2} and write
\begin{equation}
    g(x) = 
    \left(
    \begin{array}{c}
        \nabla f_1(x_1) \\
        \vdots \\
        \nabla f_N(x_N)
    \end{array}
    \right)
    =
    \left(
    \begin{array}{c}
        g_1(x_1) \\
        \vdots \\
        g_N(x_N)
    \end{array}
    \right).
    \notag
\end{equation}
\par Now we will present the framework of our proposed method, which is delineated into two primary phases: prediction and correction.
\paragraph*{Prediction step}
For a given $u^k = (x^k, \lambda^k) \in \mathcal{X} \times \Lambda$, set
\begin{subequations}
\begin{equation}
    \tilde{x}^k = P_{\mathcal{X}}[x^k - R^{-1}_k(g(x^k) - C^T \lambda^k)],
    \label{ppcm_pre_1}
\end{equation}
and
\begin{equation}
    \tilde{\lambda}^k = P_{\Lambda}[\lambda^k - S_k(C\tilde{x}^k - d)]
    \label{ppcm_pre_2}
\end{equation}
\label{ppcm_pre}
\end{subequations}
where $H_k$ is a proper chosen parameter matrix which satisfies
\begin{equation}
    \|H_k^{-1}\xi^k\|_{H_k} \leq \eta \|u^k - \tilde{u}^k\|_{H_k},
    \label{ppcm_cri}
\end{equation}
where
\begin{equation}
    \xi^k =
    \left(
    \begin{array}{c}
        g(x^k) - g(\tilde{x}^k) - C^T(\lambda^k - \tilde{\lambda}^k) \\
        0
    \end{array}
    \right).
    \label{ppcm_xi}
\end{equation}

\paragraph*{Correction step}
Calculate the new iteration $u^{k+1}$ by setting
\begin{subequations}
\begin{equation}
    x^{k+1} = P_{\mathcal{X}}[x^k - \alpha_k R_k^{-1} (g(\tilde{x}^k) - C^T \tilde{\lambda}^k)]
    \label{ppcm_cor_1}
\end{equation}
and
\begin{equation}
    \lambda^{k+1} = P_{\Lambda}[\lambda^k - \alpha_k S_k (C\tilde{x}^k - d)],
    \label{ppcm_cor_2}
\end{equation}
\label{ppcm_cor}
\end{subequations}
where
\begin{equation}
    \alpha_k = \gamma \alpha_k^*,\quad \gamma \in (0, 2) \label{ppcm_step_size},
\end{equation}
here
\begin{equation}
    \alpha_k^* = \frac{(u^k - \tilde{u}^k)^T d(u^k,\tilde{u}^k,\xi^k)}{\|d(u^k,\tilde{u}^k,\xi^k)\|^2}
    \label{ppcm_adaptive_step_size}
\end{equation}
and
\begin{equation}
    d(u^k,\tilde{u}^k,\xi^k) = H_k (u^k - \tilde{u}^k) - \xi^k. \label{ppcm_dir}
\end{equation}
\par From the above description, since both prediction and correction steps make use of projections, we call the algorithm the Projection-based Prediction-Correction Method and abbreviate it as PPCM.

\begin{remark}
In general, criterion \eqref{ppcm_cri} can be satisfied via choosing a suitable $H_k$. Since
\begin{flalign*}
    &\| H_k^{-1}\xi^k \|_{H_k}^2 \\ 
    &\qquad = \| g(x^k) - g(\tilde{x}^k) - C^T(\lambda^k - \tilde{\lambda}^k) \|_{R_k^{-1}}^2 \\
    &\qquad = \| g(x^k) - g(\tilde{x}^k) \|_{R_k^{-1}}^2 + \| C^T(\lambda^k - \tilde{\lambda}^k) \|_{R_k^{-1}}^2  \\
    &\qquad\quad - 2[g(x^k) - g(\tilde{x}^k)]^T R_k^{-1} C^T(\lambda^k - \tilde{\lambda}^k)  \\
    &\qquad \leq \| g(x^k) - g(\tilde{x}^k) \|_{R_k^{-1}}^2 + \| C^T(\lambda^k - \tilde{\lambda}^k) \|_{R_k^{-1}}^2 \\
    &\qquad\quad + 2 \| g(x^k) - g(\tilde{x}^k) \|_{R_k^{-1}} \| C^T(\lambda^k - \tilde{\lambda}^k) \|_{R_k^{-1}} \\
    &\qquad \leq \| g(x^k) - g(\tilde{x}^k) \|_{R_k^{-1}}^2 + \| C^T \|_{R_k^{-1}}^2 \| \lambda^k - \tilde{\lambda}^k \|_{R_k^{-1}}^2 \\
    &\qquad\quad + 2 \| g(x^k) - g(\tilde{x}^k) \|_{R_k^{-1}} \| C^T \|_{R_k^{-1}} \| \lambda^k - \tilde{\lambda}^k \|_{R_k^{-1}}.
\end{flalign*}
By lemma \ref{lemma_Hnorm}, we have
\begin{equation*}
    \| C^T \|_{R_k^{-1}} = \| C^T \| \leq \mu.
\end{equation*}
Synthesizing the above analysis and combining it with the fundamental inequality, we obtain
\begin{flalign*}
    &\|H_k^{-1}\xi^k\|_{H_k}^2 \\
    &\qquad \leq \| g(x^k) - g(\tilde{x}^k) \|_{R_k^{-1}}^2 + \mu^2 \|\lambda^k - \tilde{\lambda}^k\|_{R_k^{-1}}^2 \\
    &\qquad\quad + 2\mu \| g(x^k) - g(\tilde{x}^k) \|_{R_k^{-1}} \|\lambda^k - \tilde{\lambda}^k\|_{R_k^{-1}} \\
    &\qquad \leq (1 + \tau) \| g(x^k) - g(\tilde{x}^k) \|_{R_k^{-1}}^2 \\ 
    &\qquad\quad + \bigg(1 + \frac{1}{\tau}\bigg) \mu^2 \|\lambda^k - \tilde{\lambda}^k\|_{R_k^{-1}}^2, \quad \tau > 0.
\end{flalign*}
Then \eqref{ppcm_cri} can be satisfied when
\begin{equation*}
    \left\{
    \begin{aligned}
        (1 + \tau) \| g(x^k) - g(\tilde{x}^k) \|_{R_k^{-1}}^2 &\leq \eta^2 \| x^k - \tilde{x}^k \|_{R_k}^2, \\
        \bigg(1 + \frac{1}{\tau}\bigg) \mu^2 \|\lambda^k - \tilde{\lambda}^k\|_{R_k^{-1}}^2 &\leq \eta^2 \| \lambda - \tilde{\lambda}^k \|_{S_k^{-1}}^2,
    \end{aligned}
    \right.
\end{equation*}
which is equal to
\begin{equation*}
    \left\{
    \begin{aligned}
        \| g_i(x_i^k) - g_i(\tilde{x}_i^k) \| &\leq \frac{\eta r_i}{\sqrt{1 + \tau}} \| x_i^k - \tilde{x}_i^k \|, \\
        s_i &\leq \frac{\tau}{1 + \tau} \frac{\eta^2}{\mu^2} r_i,
    \end{aligned} 
    \right.
\end{equation*}
for $1 \leq i \leq N$. In fact, it is enough that $g_i$ is locally Lipschitz continuous on a neighborhood of $x_i^k$.
\end{remark}

\subsection{Theoretical analysis}
\par Consider $u^* = (x^*,\ \lambda^*) \in \Omega^*$ as an arbitrary solution to VI$(\Omega, F)$. Throughout this section, we define $u^k = (x^k,\ \lambda^k) \in \Omega$ as a given vector, $\tilde{u}^k = (\tilde{x}^k,\ \tilde{\lambda}^k) \in \Omega$ as the predictor generated by the prediction step, and $u_\alpha^{k+1} = (x_\alpha^{k+1}, \lambda_\alpha^{k+1}) \in \Omega$ as the corrector obtained from the correction step. Let
\begin{equation}
    \Theta_k(\alpha) = \|u^k - u^*\|_{H_k}^2 - \|u_\alpha^{k+1} - u^*\|_{H_k}^2
    \label{theta}
\end{equation}
which quantifies the progress achieved during the $k$-th iteration.
\par Observing that the progress $\Theta_k(\alpha)$ is dependent on the step length $\alpha$ in the correction step, it is reasonable to contemplate the maximization of this function through the selection of an optimal parameter $\alpha$. However, it is crucial to acknowledge that the solution $u^*$ remains unknown, rendering direct maximization of $\Theta_k(\alpha)$ infeasible. Consequently, the primary objective of this section is to present a lower bound for $\Theta_k(\alpha)$ that does not depend on the unknown solution $u^*$. To this end, the subsequent lemmas are dedicated.

\begin{lemma}
\par Given $u^k = (x^k,\ \lambda^k) \in \Omega$, let $\tilde{u}^k$ be the predictor obtained using \eqref{ppcm_pre}. For any $u \in \Omega$, it follows that
\begin{equation}
    (u - \tilde{u}^k)^T F(\tilde{u}^k) \geq (u - \tilde{u}^k)^T d(u^k, \tilde{u}^k, \xi^k),
    \label{eq_lem1_1}
\end{equation}
where $d(u^k, \tilde{u}^k, \xi^k)$ is defined in \eqref{ppcm_dir}.
\par Particularly, if substituting $u = u_\alpha^{k+1}$ in \eqref{eq_lem1_1}, we get
\begin{equation}
    (u_\alpha^{k+1} - \tilde{u}^k)^T F(\tilde{u}^k) \geq (u_\alpha^{k+1} - \tilde{u}^k)^T d(u^k, \tilde{u}^k, \xi^k).
    \label{eq_lem1_2}
\end{equation}
\label{lemma_1}
\end{lemma}
\begin{proof}
As $u \in \Omega$, applying Lemma \ref{lemma_proj} to \eqref{ppcm_pre}, we obtain
\begin{subequations}
    \begin{align*}
        &(x - \tilde{x}^k)^T\{\tilde{x}^k - [x^k - R_k^{-1}(g(x^k) - C^T \lambda^k)]\} \geq 0, \\
        &(\lambda - \tilde{\lambda}^k)^T\{\tilde{\lambda}^k - [\lambda^k - S_k(C\tilde{x}^k - d)]\} \geq 0.
    \end{align*}
\end{subequations}
The above can be reformulated as
\begin{subequations}
    \begin{flalign}
        &(x - \tilde{x}^k)^T[R_k(\tilde{x}^k - x^k) + (g(\tilde{x}^k) - C^T \tilde{\lambda}^k) \notag \\
        &\quad + g(x^k) - g(\tilde{x}^k) - C^T(\lambda^k - \tilde{\lambda}^k)] \geq 0, \\
        &(\lambda - \tilde{\lambda}^k)^T[S_k^{-1}(\tilde{\lambda}^k - \lambda^k) + (C\tilde{x}^k - d)] \geq 0.
    \end{flalign}
    \label{eq_lem1_3}
\end{subequations}
Using the notation of $\xi^k$, \eqref{eq_lem1_3} can be written as
\begin{equation*}
    (u - \tilde{u}^k)^T[H_k(\tilde{u}^k - u^k) + F(\tilde{u}^k) + \xi^k] \geq 0
    \label{eq_lem1_4}
\end{equation*}
and thus the assertion is proved.
\end{proof}

\begin{remark}
With the monotonicity of $F$, we have
\begin{equation*}
    (\tilde{u}^k - u^*)^T [F(\tilde{u}^k) - F(u^*)] \geq 0,
\end{equation*}
which gives that
\begin{equation*}
    (\tilde{u}^k - u^*)^T F(\tilde{u}^k) \geq (\tilde{u}^k - u^*)^T F(u^*) \geq 0.
\end{equation*}
Take $u = u^*$ in \eqref{eq_lem1_1}, we obtain
\begin{equation*}
    (u^* - \tilde{u}^k)^T F(\tilde{u}^k) \geq (u^* - \tilde{u}^k)^T d(u^k, \tilde{u}^k, \xi^k)
\end{equation*}
and thus
\begin{equation*}
    (\tilde{u}^k - u^*)^T d(u^k, \tilde{u}^k, \xi^k) \geq 0.
\end{equation*}
It follows that
\begin{flalign*}
    &(u^k - u^*)^T d(u^k, \tilde{u}^k, \xi^k) \\
    &\quad\geq\ (u^k - \tilde{u}^k)^T d(u^k, \tilde{u}^k, \xi^k) \\
    &\quad\stackrel{\mathclap{\eqref{ppcm_dir}}}{=}\ \|u^k - \tilde{u}^k\|_{H_k}^2 - (u^k - \tilde{u}^k)^T \xi^k \\
    &\quad\geq\ \|u^k - \tilde{u}^k\|_{H_k}^2 - \|u^k - \tilde{u}^k\|_{H_k} \|H_k^{-1}\xi^k\|_{H_k} \\
    &\quad\stackrel{\mathclap{\eqref{ppcm_cri}}}{\geq}\ (1 - \eta) \|u^k - \tilde{u}^k\|_{H_k}^2 \\
    &\quad>\ 0. \text{\quad (when $u^k \neq \tilde{u}^k$)}
\end{flalign*}
Therefore, the quantity $d(u^k, \tilde{u}^k, \xi^k)$ defined in \eqref{ppcm_dir} represents an ascending direction of the unknown function $\frac{1}{2}\|u - u^*\|^2$ under Euclidean norm at point $u^k$, and hence $H_k^{-1}d(u^k, \tilde{u}^k, \xi^k)$ gives an ascending direction of $\frac{1}{2}\|u - u^*\|_{H_k}^2$. This greatly aids us in proving the convergence of the algorithm.
\end{remark}

\begin{lemma}
Given $u^k = (x^k, \lambda^k) \in \Omega$, let $\tilde{u}^k$ be the predictor obtained using \eqref{ppcm_pre}, and $u_\alpha^{k+1}$ be the corrector (dependent on $\alpha$) produced by \eqref{ppcm_cor}. Then we have
\begin{equation}
    \Theta_k(\alpha) \geq \|u_\alpha^{k+1} - u^k\|_{H_k}^2 + 2 \alpha (u_\alpha^{k+1} - \tilde{u}^k)^T F(\tilde{u}^k).
    \label{eq_lem2_1}
\end{equation}
\label{lemma_2}
\end{lemma}
\begin{proof}
Since $u^* \in \Omega$, applying Lemma \ref{lemma_proj} to \eqref{ppcm_cor} yields
\begin{subequations}
    {\small
    \begin{flalign}
        &(x^* - x_\alpha^{k+1})^T\{x_\alpha^{k+1} - [x^k - \alpha R_k^{-1}(g(\tilde{x}^k) - C^T \tilde{\lambda}^k)]\} \geq 0, \\
        &(\lambda^* - \lambda_\alpha^{k+1})^T\{\lambda_\alpha^{k+1} - [\lambda^k - \alpha S_k(C\tilde{x}^k - d)]\} \geq 0.
    \end{flalign}
    }
    \label{eq_lem2_2}
\end{subequations}
From \eqref{eq_lem2_2}, we obtain
\begin{equation*}
    (u^* - u_\alpha^{k+1})^T[H_k(u_\alpha^{k+1} - u^k) + \alpha F(\tilde{u}^k)] \geq 0
\end{equation*}
and thus
\begin{equation*}
    (u^* - u_\alpha^{k+1})^T H_k(u_\alpha^{k+1} - u^k) \geq \alpha(u_\alpha^{k+1} - u^*)^T F(\tilde{u}^k).
\end{equation*}
Notice that we have the following identity:
\begin{flalign*}
    &2(u^* - u_\alpha^{k+1})^T H_k(u_\alpha^{k+1} - u^k) \\
    &\quad = \|u^* - u^k\|_{H_k}^2 - \|u^* - u_\alpha^{k+1}\|_{H_k}^2 - \|u_\alpha^{k+1} - u^k\|_{H_k}^2 \\
    &\quad = \Theta_k(\alpha) - \|u_\alpha^{k+1} - u^k\|_{H_k}^2.
\end{flalign*}
Therefore,
{\footnotesize
\begin{flalign*}
    \Theta_k(\alpha) &= \|u_\alpha^{k+1} - u^k\|_{H_k}^2 + 2(u^* - u_\alpha^{k+1})^T H_k(u_\alpha^{k+1} - u^k) \\
    &\geq \|u_\alpha^{k+1} - u^k\|_{H_k}^2 + 2\alpha(u_\alpha^{k+1} - u^*)^T F(\tilde{u}^k) \\
    &= \|u_\alpha^{k+1} - u^k\|_{H_k}^2 + 2\alpha[(u_\alpha^{k+1} - \tilde{u}^k) + (\tilde{u}^k - u^*)]^T F(\tilde{u}^k) \\
    &\geq \|u_\alpha^{k+1} - u^k\|_{H_k}^2 + 2\alpha(u_\alpha^{k+1} - \tilde{u}^k)^T F(\tilde{u}^k).
\end{flalign*}
}
The last inequality follows from the monotonicity of $F$ and $(\tilde{u}^k - u^*)^T F(u^*) \geq 0$. Consequently, the assertion of this lemma is readily derived from the aforementioned inequality.
\end{proof}

\par Hence, we arrive at the following lemma.
\begin{lemma}
Given $u^k = (x^k, \lambda^k) \in \Omega$, let $\tilde{u}^k$ be the predictor obtained using \eqref{ppcm_pre}, and $u_\alpha^{k+1}$ be the corrector (dependent on $\alpha$) produced by \eqref{ppcm_cor}. Then we have
{\small
\begin{equation}
    \Theta_k(\alpha) \geq \|u_\alpha^{k+1} - u^k\|_{H_k}^2 + 2 \alpha (u_\alpha^{k+1} - \tilde{u}^k)^T d(u^k, \tilde{u}^k, \xi^k).
    \label{eq_thm_theta}
\end{equation}
}
\label{thm_theta}
\end{lemma}
\begin{proof}
The assertion follows from Lemma \ref{lemma_1} and Lemma \ref{lemma_2} directly.
\end{proof}

\par The subsequent lemma furnishes a concave quadratic function of $\alpha$ that serves as a lower bound for $\Theta_k(\alpha)$.
\begin{lemma}
Let $\Theta(\alpha)$ be defined by \eqref{theta} and $d(u^k, \tilde{u}^k, \xi^k)$ be defined by \eqref{ppcm_dir}. Then for any $u^* \in \Omega^*$ and $\alpha > 0$, we have
\begin{equation}
    \Theta_k(\alpha) \geq \Phi_k(\alpha),
    \label{eq_thm_phi_1}
\end{equation}
where
{\small
\begin{equation}
    \Phi_k(\alpha) := 2\alpha(u^k - \tilde{u}^k)^T d(u^k, \tilde{u}^k, \xi^k) - \alpha^2 \|H_k^{-1}d(u^k, \tilde{u}^k, \xi^k)\|_{H_k}^2.
    \label{eq_thm_phi_2}
\end{equation}
}
\label{thm_phi}
\end{lemma}
\begin{proof}
It follows from (\ref{eq_thm_theta}) that
\begin{flalign*}
    \Theta_k(\alpha) &\geq  2 \alpha [(u_\alpha^{k+1} - u^k) + (u^k - \tilde{u}^k)]^T d(u^k, \tilde{u}^k, \xi^k) \\ 
    &\quad + \|u^k - u_\alpha^{k+1}\|_{H_k}^2 \\
    &= 2 \alpha (u^k - \tilde{u}^k)^T d(u^k, \tilde{u}^k, \xi^k) \\
    &\quad + \|u^k - u_\alpha^{k+1}\|_{H_k}^2 - 2 \alpha (u^k - u_\alpha^{k+1})^T d(u^k, \tilde{u}^k, \xi^k) \\
    &= 2 \alpha (u^k - \tilde{u}^k)^T d(u^k, \tilde{u}^k, \xi^k) \\
    &\quad - \alpha^2 \|H_k^{-1}d(u^k, \tilde{u}^k, \xi^k)\|_{H_k}^2 \\
    &\quad + \|u^k - u_\alpha^{k+1} - \alpha H_k^{-1} d(u^k, \tilde{u}^k, \xi^k)\|_{H_k}^2 \\
    &\geq 2\alpha(u^k - \tilde{u}^k)^T d(u^k, \tilde{u}^k, \xi^k) \\
    &\quad - \alpha^2 \|H_k^{-1}d(u^k, \tilde{u}^k, \xi^k)\|_{H_k}^2.
\end{flalign*}
The assertion follows from the definition of $\Phi_k(\alpha)$ directly.
\end{proof}

\par Observe that $\Theta_k(\alpha)$ can be considered the progress achieved by $u_\alpha^{k+1}$, and $\Phi_k(\alpha)$ serves as a lower bound for $\Theta_k(\alpha)$. Consequently, it inspires us to determine an appropriate $\alpha$ that maximizes $\Phi_k(\alpha)$. As $\Phi_k(\alpha)$ is a concave quadratic function of $\alpha$, its maximum occurs at
\begin{equation}
    \alpha_k^* = \frac{(u^k - \tilde{u}^k)^T d(u^k, \tilde{u}^k, \xi^k)}{\|H_k^{-1}d(u^k, \tilde{u}^k, \xi^k)\|_{H_k}^2}.
    \label{alpha_star}
\end{equation}
This corresponds precisely to the same $\alpha_k^*$ as in \eqref{ppcm_adaptive_step_size}. Note that
\begin{equation}
    \Phi_k(\alpha_k^*) = \alpha_k^*(u^k - \tilde{u}^k)^T d(u^k, \tilde{u}^k, \xi^k).
    \label{phi_alpha_star}
\end{equation}
\par It is pertinent to mention that, under the condition stated in \eqref{ppcm_cri}, we have
\begin{flalign}
    &2(u^k - \tilde{u}^k)^T d(u^k, \tilde{u}^k, \xi^k) - \|H_k^{-1}d(u^k, \tilde{u}^k, \xi^k)\|_{H_k}^2 \notag \\
    &\quad = \|u^k - \tilde{u}^k\|_{H_k}^2 - \|u^k - \tilde{u}^k - H_k^{-1}d(u^k, \tilde{u}^k, \xi^k)\|_{H_k}^2  \notag \\
    &\quad = \|u^k - \tilde{u}^k\|_{H_k}^2 - \|H_k^{-1}\xi^k\|_{H_k}^2 \notag \\
    &\quad \geq (1 - \eta^2)\|u^k - \tilde{u}^k\|_{H_k}^2. \label{eq_4.26}
\end{flalign}
Therefore, whenever $u^k \neq \tilde{u}^k$, it follows from \eqref{alpha_star} and \eqref{eq_4.26} that
\begin{equation}
    \alpha_k^* > \frac{1}{2}.
    \label{eq_alpha_star}
\end{equation}
Consequently, from \eqref{alpha_star}, \eqref{phi_alpha_star} and \eqref{eq_alpha_star} we obtain
\begin{equation}
    \Phi_k(\alpha_k^*) \geq \frac{1 - \eta^2}{4} \|u^k - \tilde{u}^k\|_{H_k}^2.
    \label{eq_phi_alpha_star}
\end{equation}
In light of numerical experiments, we find it advantageous to scale the ``optimal" value $\alpha_k^*$ by a relaxation factor $\gamma \in (0, 2)$ (preferably lies in $[1, 2)$). As such, we recommend employing the correction formula \eqref{ppcm_cor} with the step-size defined in \eqref{ppcm_step_size}.
\begin{remark}
In a distributed setting, utilizing the step size calculated by \eqref{ppcm_step_size} would preclude decentralization. However, it is noteworthy that there always exists $\gamma_k \in (0, 2)$ such that $\gamma_k \alpha_k^* = 1$, allowing for the direct use of a unit step size $\alpha_k = 1$. This adaptation ensures the algorithm's seamless application in distributed environments.
\end{remark}

\par Now, our primary focus lies in investigating the convergence of the proposed method. The ensuing lemma addresses the contractive property of the sequence generated by the proposed method.
\begin{lemma}
Let $u^{k+1} = u^{k+1}(\gamma \alpha_k^*)$ be the new iteration. Then for any $u^* \in \Omega^*$ and $\gamma \in (0, 2)$, we have
{\small
\begin{equation}
    \|u^{k+1} - u^*\|_{H_k}^2 \leq \|u^k - u^*\|_{H_k}^2 - \frac{\gamma(2 - \gamma)(1 - \eta^2)}{4}\|u^k - \tilde{u}^k\|_{H_k}^2.
    \label{eq_thm_con_1}
\end{equation}
}
\label{thm_con}
\end{lemma}
\begin{proof}
Through a straightforward manipulation, we obtain
\begin{flalign}
    \Phi_k(\gamma\alpha_k^*) &= 2\gamma\alpha_k^*(u^k - \tilde{u}^k)^T d(u^k, \tilde{u}^k, \xi^k) \notag \\
    &\quad - (\gamma^2\alpha_k^*) (\alpha_k^* \|H_k^{-1}d(u^k, \tilde{u}^k, \xi^k)\|_{H_k}^2) \notag \\
    &= (2 \gamma \alpha_k^* - \gamma^2 \alpha_k^*) (u^k - \tilde{u}^k)^T d(u^k, \tilde{u}^k, \xi^k) \notag \\
    &= \gamma(2 - \gamma) \Phi_k(\alpha_k^*). \label{eq_thm_con_2}
\end{flalign}
It follows from Lemma \ref{thm_phi}, \eqref{eq_phi_alpha_star} and \eqref{eq_thm_con_2} that
\begin{flalign*}
    \Theta_k(\gamma\alpha_k^*) &= \|u^k - u^*\|_{H_k}^2 - \|u^{k+1}(\gamma\alpha_k^*) - u^*\|_{H_k}^2 \\ 
    &\geq \Phi_k(\gamma\alpha_k^*) \\
    &= \gamma(2 - \gamma) \Phi_k(\alpha_k^*) \\
    &\geq \frac{\gamma(2 - \gamma)(1 - \eta^2)}{4}\|u^k - \tilde{u}^k\|_{H_k}^2.
\end{flalign*}
\end{proof}

\par According to Lemma \ref{thm_con}, there exists a constant $c > 0$ such that
\begin{equation*}
    \|u^{k+1} - u^*\|_{H_k}^2 \leq \|u^k - u^*\|_{H_k}^2 - c\|u^k - \tilde{u}^k\|_{H_k}^2, \quad \forall u^* \in \Omega^*.
\end{equation*}
As a result, the sequence ${u^k}$ is bounded, and the proposed method is categorized as a contractive method since its new iteration is closer to the solution set $\Omega^*$. With this understanding, we proceed to demonstrate the convergence of the proposed method.

\begin{theorem}
The sequence $\{u^k\}$ generated by the proposed method converges to some $u^{\infty}$ which is a solution of VI$(\Omega, F)$.
\label{thm_cvg}
\end{theorem}
\begin{proof}
By virtue of Lemma \ref{thm_con}, the sequence $\{u^k\}$ is bounded, and
\begin{equation*}
    \lim\limits_{k \rightarrow \infty} \|u^k - \tilde{u}^k\|_{H_k} = 0.
\end{equation*}
Consequently, $\{\tilde{u}^k\}$ is also bounded. Combining the above with \eqref{ppcm_cri}, we can deduce that
\begin{equation*}
    \lim\limits_{k \rightarrow \infty} d(u^k, \tilde{u}^k, \xi^k) = 0.
\end{equation*}
As a result, employing Lemma \ref{lemma_1}, we obtain that for all $u \in \Omega$,
\begin{equation*}
    \lim\limits_{k \rightarrow \infty} (u - \tilde{u}^k)^T F(\tilde{u}^k) \geq \lim\limits_{k \rightarrow \infty} (u - \tilde{u}^k)^T d(u^k, \tilde{u}^k, \xi^k) = 0.
\end{equation*}
As $\{\tilde{u}^k\}$ is bounded, it possesses at least one cluster point. Let $u^\infty$ denote a cluster point of $\{\tilde{u}^k\}$, and the subsequence $\{\tilde{u}^{k_j}\}$ converges to $u^\infty$. Therefore, we have
\begin{equation*}
    \lim\limits_{j \rightarrow \infty} (u - \tilde{u}^{k_j})^T F(\tilde{u}^{k_j}) \geq 0,\quad \forall u \in \Omega
\end{equation*}
and consequently
\begin{equation*}
    (u - u^\infty)^T F(u^\infty) \geq 0,\quad \forall u \in \Omega.
\end{equation*}
This means that $u^\infty$ is a solution of VI$(\Omega, F)$. Substituting $u^*$ by $u^\infty$ in \eqref{eq_thm_con_1}, we obtain
\begin{equation}
    \|u^{k+1} - u^\infty\|_{H_k}^2 \leq \|u^k - u^\infty\|_{H_k}^2,\quad \forall k \geq 0.
    \label{eq_thm_cvg_1}
\end{equation}
Since $\tilde{u}^{k_j} \rightarrow u^\infty$ as $j \rightarrow \infty$ and $u^k - \tilde{u}^k \rightarrow 0$ as $k \rightarrow \infty$, then for any given $\epsilon > 0$, there exists an integer $l > 0$, such that
\begin{equation}
    \|\tilde{u}^{k_l} - u^\infty\|_{H_k} < \frac{\epsilon}{2}\quad \text{and} \quad \|u^{k_l} - \tilde{u}^{k_l}\|_{H_k} < \frac{\epsilon}{2}.
    \label{eq_thm_cvg_2}
\end{equation}
Therefore, for any $k \geq k_l$, it follows from \eqref{eq_thm_cvg_1} and \eqref{eq_thm_cvg_2} that
\begin{align*}
    \|u^k - u^\infty\|_{H_k} &\leq \|u^{k_l} - u^\infty\|_{H_k} \\
    &\leq \|u^{k_l} - \tilde{u}^{k_l}\|_{H_k} + \|\tilde{u}^{k_l} - u^\infty\|_{H_k} \\
    &< \epsilon.
\end{align*}
This implies that the sequence $\{u^k\}$ converges to $u^\infty$ which is a solution of VI$(\Omega, F)$.
\end{proof}

\section{PPCM for Distributed Consensus Optimization}
\subsection{Algorithm in detail}
\par Now, we will elaborate on the specific steps of the Projection-based Prediction-Correction Method (PPCM) for (P). As mentioned in the previous section, we will adopt a unit step size $\alpha_k = 1$ to enhance the algorithm's adaptability to distributed environments. Recall that $A = (a_{ij})_{N \times N}$ represents the adjacency matrix of the undirected graph corresponding to the communication network. Let $\eta \in (0, 1)$, $\tau > 0$ and $\mu$ be given constants known to all agents, with $2\| A \|_{\infty} \leq \mu$ (which can be achieved by adjusting $A$). Since $C = L \otimes I_n$, we have $\| C^T \| = \| L^T \| = \rho(L) \leq 2 \| A \|_{\infty} \leq \mu$. The initial state $\lambda_i^0$ is set to the zero vector.

\paragraph*{Prediction step}
For node $i$, calculate
\begin{equation*}
    \tilde{x}_i^k = P_{\mathcal{X}_i}\bigg[x_i^k - \frac{1}{r_i^k}\bigg(g_i(x_i^k) - \sum\limits_{j \in \mathcal{N}_i}a_{ij}(\lambda_i^k - \lambda_j^k)\bigg)\bigg],
\end{equation*}
where $r_i^k$ is a suitably chosen parameter that satisfies
\begin{equation*}
    \|g(x_i^k) - g(\tilde{x}_i^k)\| \leq \frac{\eta r_i^k}{\sqrt{1 + \tau}} \|x_i^k - \tilde{x}_i^k\|.
\end{equation*}
Subsequently, communicate with neighboring nodes to exchange information, obtaining all $\tilde{x}_j^k$ for $j \in \mathcal{N}_i$, and compute
\begin{equation*}
    \lambda_i^{k+1} = \lambda_i^k - \frac{\eta^2 \tau}{\mu^2 (1 + \tau)} r_i^k \sum\limits_{j \in \mathcal{N}_i}a_{ij}(\tilde{x}_i^k - \tilde{x}_j^k).
\end{equation*}

\paragraph*{Correction step}
For node $i$, communicate again with adjacent nodes to exchange information, obtaining all $\lambda_j^{k+1}$ for $j \in \mathcal{N}_i$. Subsequently, compute
\begin{equation}
    x_i^{k+1} = P_{\mathcal{X}_i}\bigg[x_i^k - \frac{1}{r_i^k}\bigg(g_i(\tilde{x}_i^k) - \sum\limits_{j \in \mathcal{N}_i}a_{ij}(\lambda_i^{k+1} - \lambda_j^{k+1})\bigg)\bigg]. \notag
\end{equation}

\begin{remark}
Note that the selection of the adjacency matrix $A = (a_{ij})_{p \times p}$ is not unique. One viable option is:
\begin{equation*}
    a_{ij} = \frac{\tau}{2(1 + \tau)} *
    \begin{cases}
        \frac{1}{\max\{d_i, d_j\}} & (i, j) \in E, \\
        0 & (i, j) \notin E.
    \end{cases}
\end{equation*}
where $d_i = |\mathcal{N}_i|$ is the degree of node $i$. It is not difficult to derive that $\| A \|_{\infty} \leq \frac{\tau}{2(1 + \tau)} \leq \frac{1}{2}\sqrt{\frac{\tau}{1 + \tau}}$. Therefore, one can take $\mu = \sqrt{\frac{\tau}{1 + \tau}}$.
\end{remark}

\begin{remark}
The choice of $\mu$ here is also not unique. The selection adopted in this paper aims to satisfy the condition $\frac{\tau}{\mu^2 (1 + \tau)} = 1$, thereby simplifying the step of computing $\lambda_i^{k+1}$ and reducing the number of uniform parameters that nodes need to preset.
\end{remark}

\par Besides, in practical calculations, the utilization of a self-adjustment technique for parameter tuning becomes pivotal, exerting a profound influence on computational efficiency. We hereby provide a more specific pseudocode of PPCM.

\begin{algorithm}
    \caption{Distributed Consensus PPCM for agent $i$}
    \label{algo:DC-PPCM}
    \Initialize $k = 0,\ \mathrm{tol} = 10^{-3},\ \epsilon = 1.0,\ \eta = 0.9,\ \tau = 1.5,\ r_i^0 = 1.0,\ x_i^0 \in \mathcal{X}_i,\ \lambda_i^0 = 0$\;
    compute $a_{ij} = \frac{\tau}{2(1+\tau)} \frac{1}{\max\{d_i, d_j\}}$ for $j \in \mathcal{N}_i$\; 
    \While{$\epsilon < \mathrm{tol}$}{
      \SetNoFillComment
      \Comment{Adaptive parameter tuning}
      \While{$\mathrm{True}$}{
        {\small
        $\tilde{x}_i^k = P_{\mathcal{X}_i}\bigg[x_i^k - \frac{1}{r_i^k}\bigg(g_i(x_i^k) - \sum\limits_{j \in \mathcal{N}_i} a_{ij}(\lambda_i^k - \lambda_j^k) \bigg)\bigg]$\;
        }
        $t = \frac{\sqrt{1 + \tau} \|g(x_i^k) - g(\tilde{x}_i^k)\|}{r_i^k\|x_i^k - \tilde{x}_i^k\|}$\;
        \eIf{$t > \eta$}{
          $r_i^k = r_i^k * 1.5 * \max\{1, t\}$\;
        }{
          \Break\;
        }
      }
      communicate with neighboring nodes to obtain $\tilde{x}_j^k$ for $j \in \mathcal{N}_i$\;
      $s_i^k = \eta^2 r_i^k$\;
      $\lambda_i^{k+1} = \lambda_i^k - s_i^k \sum\limits_{j \in \mathcal{N}_i} a_{ij}(\tilde{x}_i^k - \tilde{x}_j^k)$\;
      communicate with neighboring nodes to obtain $\lambda_j^{k+1}$ for $j \in \mathcal{N}_i$\;
      {\footnotesize
      $x_i^{k+1} = P_{\mathcal{X}_i}\bigg[x_i^k - \frac{1}{r_i^k}\bigg(g_i(\tilde{x}_i^k) - \sum\limits_{j \in \mathcal{N}_i} a_{ij}(\lambda_i^{k+1} - \lambda_j^{k+1})\bigg)\bigg]$\;
      }
      $\epsilon = \max\{\sqrt{r_i^k}\|x_i^k - \tilde{x}_i^k\|_{\infty},\ \|\lambda_i^k - \lambda_i^{k+1}\|_{\infty} / \sqrt{s_i^k} \}$\;
      \eIf{$t \leq 0.5$}{
        $r_i^{k+1} = r_i^k * \frac{t}{0.7}$\;
      }{
        $r_i^{k+1} = r_i^k$\;
      }
      $k = k + 1$\;
    }
    \KwRet{$x_i^k$}
\end{algorithm}

\subsection{Numerical experiments}
\par To verify the correctness of the PPCM and to delve into its performance among other aspects, we will explain these issues in detail through three specific examples: linear least squares, logistic regression, and support vector machines. We simulate a distributed environment locally using Python's mpi4py library, with communication networks divided into ring and fully connected (P2P) types, corresponding to ring graphs and complete graphs in graph theory, respectively. Our codes are publicly available at \url{https://github.com/harmoke/Distributed-Optimization}, and all code runs on a Macbook Pro M3 Max.
\par At the same time, we have also implemented the algorithm proposed by Nedi\'{c} et al. \cite{nedic2010constrained}, which has provided us with valuable references. Given that this algorithm belongs to the gradient methods utilizing weighted averages, for the sake of concise expression, we refer to it as the Weighted Averaging Gradient Method (WAGM). Its specific form is as follows:
\begin{equation}
    \left\{
    \begin{aligned}
        &\tilde{x}_i^k = \sum\limits_{j \in \mathcal{N}_i} w_{ij} x_i^k, \\
        &x_i^{k+1} = P_{\mathcal{X}_i}[\tilde{x}_i^k - \alpha_k g(\tilde{x}_i^k)],
    \end{aligned}
    \right.
    \label{method_dgm}
\end{equation}
here $W = (w_{ij})_{N \times N}$ denotes a positive doubly stochastic matrix (the sum of each row and column equals 1), which matches the structure of graph $\mathcal{G}$. Furthermore, the step size should satisfy the following conditions: for all $k \geq 1$, $0 < \alpha_{k+1} \leq \alpha_k$; $\sum\limits_{k=1}^\infty \alpha_k = \infty$; and $\sum\limits_{k=1}^\infty \alpha_k^2 < \infty$. For more theoretical details, refer to \cite{nedic2018distributed}. 
\begin{remark}
The choice of $W$ is not unique. We adopted a most straightforward method:
\begin{equation*}
    w_{ij} = \frac{1}{N} *
    \begin{cases}
        1 & (i, j) \in E, \\
        N - d_i & i = j, \\
        0 & otherwise.
    \end{cases}
\end{equation*}
The step size $\alpha_k$ plays a critical role in the convergence process of WAGM, but theoretically, no clear guidance is provided, thus it must be manually adjusted according to specific computation problems. Through a series of detailed experiments, we found relatively appropriate step size parameter settings, which will be presented in the following sections along with specific experiments.
\end{remark}
We now provide more specific pseudocode for WAGM:
\begin{algorithm}
    \caption{Distributed Consensus WAGM for agent $i$}
    \label{algo:DC-WAGM}
    \Initialize $k = 0,\ \mathrm{tol} = 10^{-6},\ \epsilon = 1.0,\ x_i^0 \in \mathcal{X}_i$\;
    \While{$\epsilon < \mathrm{tol}$}{
      communicate with neighboring nodes to obtain $x_j^k$ for $j \in \mathcal{N}_i$\;
      $\tilde{x}_i^k = \Big(1 - \frac{d_i}{N}\Big)x_i^k + \frac{1}{N}\sum\limits_{j \in \mathcal{N}_i} x_j^k$\;
      $x_i^{k+1} = P_{\mathcal{X}_i}\Big[\tilde{x}_i^k - \alpha_k g_i(\tilde{x}_i^k)\Big]$\;
      $\epsilon = \| x_i^{k+1} - x_i^k \|_{\infty}$\;
      $k = k + 1$\;
    }
    \KwRet{$x_i^k$}
\end{algorithm}

\subsubsection{Linear Least Squares}
\par Solving for the least squares solution of a system of linear equations is a classic and fundamental problem, for which many mature methods have been developed. However, when the scale of the equation system significantly increases or when the solution needs to be found within a decentralized distributed computing framework, this seemingly simple problem immediately becomes a task of great challenge and complexity.
\par Given $Q \in \mathbb{R}^{m \times n}$ and $y \in \mathbb{R}^{m \times 1}$, finding the least squares solution of the linear system $Qx = y$ essentially involves solving the following optimization problem:
\begin{equation*}
   \min\limits_{x \in \mathbb{R}^n}\ \frac{1}{2}\|Qx - y\|^2.
\end{equation*}
Next, we divide $Q$ and $y$ into $N$ vertical blocks, namely
\begin{equation*}
    Q = \left(
    \begin{array}{c}
        Q_1 \\
        \vdots \\
        Q_N
    \end{array}
    \right) \in \mathbb{R}^{m \times n}, \qquad
    y = \left(
    \begin{array}{c}
        y_1 \\
        \vdots \\
        y_N
    \end{array}
    \right) \in \mathbb{R}^{m \times 1},
\end{equation*}
where $Q_i \in \mathbb{R}^{m_i \times n}$ and $y_i \in \mathbb{R}^{m_i \times 1}$, and $m = \sum\limits_{i = 1}^N m_i$. Note that
{\small
\begin{equation*}
    \min\limits_{x \in \mathbb{R}^n}\ \frac{1}{2}\|Qx - y\|^2 =  \min\limits_{x \in \mathbb{R}^n}\ \sum\limits_{i=1}^N \frac{1}{2}\|Q_i x - y_i\|^2 = \min\limits_{x \in \mathbb{R}^n}\ \sum\limits_{i=1}^N f_i(x),
\end{equation*}
}
where
\begin{equation*}
    f_i(x) = \frac{1}{2}\|Q_i x - y_i\|^2.
\end{equation*}
Therefore, by assigning the matrix-vector pair $(Q_i, y_i)$ to agent $i$, we are able to solve the problem using distributed consensus optimization algorithms.

\paragraph*{Experimental Setup}
Set $n = 63000,\ m = 4000$, and randomly select the matrix $Q$ and the vector $y$ from the standard normal distribution. The communication network consists of $N$ nodes in either a ring or a complete graph ($2 \leq N \leq 10$). The standard solution is obtained using the lstsq function from the Numpy library.

\paragraph*{PPCM Parameter Selection}
Parameters are chosen according to the default values provided in Algorithm \ref{algo:DC-PPCM}.

\paragraph*{WAGM Parameter Selection}
Set $\alpha_k = \frac{10^{-4}}{k + 1}$, with other parameters chosen according to the default values provided in Algorithm \ref{algo:DC-WAGM}.

\subsubsection{Logistic Regression}
\par Logistic regression is a widely used method in the fields of machine learning and statistics, with its core problem being the solution of the following optimization problem:
\begin{equation*}
    \max\limits_{x \in \mathbb{R}^n}\ \frac{1}{M}\sum\limits_{i=1}^M [b_i * \log(\sigma(a_i^T x)) + (1 - b_i) * \log(1 - \sigma(a_i^T x))],
\end{equation*}
where $a_i \in \mathbb{R}^n$ represents data points, $b_i \in \{0, 1\}$ are the corresponding labels, and
\begin{equation*}
    \sigma(z) = \frac{1}{1 + \exp(-z)}
\end{equation*}
is the Sigmoid function. For distributed solving, the $M$ data points are divided into $N$ groups.

\paragraph*{Experimental Setup}
Set $M = 5000$, and generate $(a_i,\ b_i)_{i=1}^M$ using the make\_classification function from the Sklearn library. The communication network is set to a complete graph with $N=5$. The standard solution is obtained using the minimize function from the Scipy library.

\paragraph*{PPCM Parameter Selection}
Set $\mathrm{tol} = 10^{-6}$, with other parameters chosen according to the default values provided in Algorithm \ref{algo:DC-PPCM}.

\paragraph*{WAGM Parameter Selection}
Set $\alpha_k = \frac{20}{k + 1}$, with other parameters chosen according to the default values provided in Algorithm \ref{algo:DC-WAGM}.

\subsubsection{Support Vector Machine}
\par Support Vector Machine (SVM) is a powerful supervised learning model widely applied in classification and regression problems. The core idea of SVM is to find the optimal hyperplane in the feature space that separates data points of different categories. This optimal hyperplane is defined as the one that maximizes the margin between categories of data points, thereby providing a clear decision boundary.
\par Considering the linear soft-margin SVM for binary classification problems, its mathematical model can be expressed as follows:
\begin{equation*}
    \min\limits_{x \in \mathbb{R}^n}\ \frac{1}{2}\|x\|^2 + \frac{\theta}{M} \sum\limits_{i=1}^M \max\{0, 1 - b_i a_i^T x\},
 \end{equation*}
 where $a_i \in \mathbb{R}^n$ represents data points, $b_i \in \{-1, 1\}$ are the corresponding labels, and $\theta$ is a weight parameter. For distributed solving, the $M$ data points are divided into $N$ groups.

 \paragraph*{Experimental Setup}
Set $M = 10000$, and generate $(a_i, b_i)_{i=1}^M$ using the make\_classification function from the Sklearn library. Additionally, a convex set constraint $x \geq 0$ is added. The communication network is set to a complete graph with $N=5$. The standard solution is obtained using the minimize function from the Scipy library.

\paragraph*{PPCM Parameter Selection}
Set $\mathrm{tol} = 10^{-5}$, with other parameters chosen according to the default values provided in Algorithm \ref{algo:DC-PPCM}.

\paragraph*{WAGM Parameter Selection}
Set $\alpha_k = \frac{20}{k + 1}$, with other parameters chosen according to the default values provided in Algorithm \ref{algo:DC-WAGM}.

\subsection{Experimental Results and Analysis}
\par Upon a detailed analysis of the data presented in several tables, it is evident that PPCM demonstrates a significant efficiency advantage over WAGM and built-in library functions. Specifically, PPCM not only achieves a speed improvement of at least four times compared to WAGM but also exhibits an acceleration ratio up to twelve times when compared to built-in library functions. Notably, while significantly accelerating computation speed, PPCM still maintains a high level of accuracy, whereas WAGM suffers a substantial loss in precision. Despite PPCM requiring two rounds of communication per iteration, compared to the single round needed by WAGM, the total communication cost of PPCM is actually lower due to its fewer required iterations.

\par Further examination of Tables \ref{tab:leastsquare_cycle} and \ref{tab:leastsquare_complete} reveals that as the number of nodes increases, so does the required runtime. This phenomenon is primarily due to the expanded network size, which leads to higher communication overheads.

\par A comparison of the data from these two tables further reveals that, compared to the ring graph, the complete graph structure shows significant improvements in efficiency and accuracy. This is because, in a more tightly knit network graph, each node receives a larger and more comprehensive amount of information in each communication round, greatly facilitating rapid convergence of the algorithm.

\par Moreover, Tables \ref{tab:logistic} and \ref{tab:svm} showcase PPCM's exceptional performance in distributed logistic regression and support vector machine problems.

\par Overall, PPCM, with its high efficiency, accuracy, flexibility, and robustness, proves to be an exceptionally superior method.

\begin{table*}[]
    \centering
    \caption{Distributed Linear Least Squares (Ring)}
    \label{tab:leastsquare_cycle}
    \resizebox{0.9\textwidth}{!}{%
    \begin{tabular}{|c|c|c|c|cc|}
    \hline
    \multirow{2}{*}{Methods}   & \multirow{2}{*}{N} & \multirow{2}{*}{Average steps} & \multirow{2}{*}{Average running time/s} & \multicolumn{2}{c|}{Average error}                                                \\ \cline{5-6} 
                          &                      &                         &                               & \multicolumn{1}{c|}{$L_2$ norm}                  & $L_\infty$ norm             \\ \hline
    \multirow{9}{*}{PPCM} & 2                    & $30$                    & $3.00541$                     & \multicolumn{1}{c|}{$7.50367 \times 10^{-4}$} & $4.32202 \times 10^{-5}$ \\ \cline{2-6} 
                          & 3                    & $28$                    & $6.43045$                     & \multicolumn{1}{c|}{$7.13905 \times 10^{-4}$} & $4.45913 \times 10^{-5}$ \\ \cline{2-6} 
                          & 4                    & $44$                    & $13.0285$                     & \multicolumn{1}{c|}{$1.65427 \times 10^{-4}$} & $9.68939 \times 10^{-6}$ \\ \cline{2-6} 
                          & 5                    & $52$                    & $19.7398$                     & \multicolumn{1}{c|}{$4.03915 \times 10^{-4}$} & $2.37027 \times 10^{-5}$ \\ \cline{2-6} 
                          & 6                    & $78$                    & $32.9723$                     & \multicolumn{1}{c|}{$4.96195 \times 10^{-4}$} & $2.93456 \times 10^{-5}$ \\ \cline{2-6} 
                          & 7                    & $123$                   & $67.7719$                     & \multicolumn{1}{c|}{$6.98445 \times 10^{-4}$} & $4.03916 \times 10^{-5}$ \\ \cline{2-6} 
                          & 8                    & $208$                   & $102.996$                     & \multicolumn{1}{c|}{$1.14660 \times 10^{-3}$} & $6.75993 \times 10^{-5}$ \\ \cline{2-6} 
                          & 9                    & $309$                   & $175.169$                     & \multicolumn{1}{c|}{$1.43015 \times 10^{-3}$} & $8.50363 \times 10^{-5}$ \\ \cline{2-6} 
                          & 10                   & $425$                   & $251.392$                     & \multicolumn{1}{c|}{$1.96020 \times 10^{-3}$} & $1.12042 \times 10^{-4}$ \\ \hline
    \multirow{9}{*}{WAGM} & 2                    & $216$                   & $13.2036$                     & \multicolumn{1}{c|}{$3.61078 \times 10^{-3}$} & $2.13656 \times 10^{-4}$ \\ \cline{2-6} 
                          & 3                    & $249$                   & $20.8617$                     & \multicolumn{1}{c|}{$2.77573 \times 10^{-3}$} & $1.61357 \times 10^{-4}$ \\ \cline{2-6} 
                          & 4                    & $307$                   & $43.0887$                     & \multicolumn{1}{c|}{$4.13048 \times 10^{-3}$} & $2.58122 \times 10^{-4}$ \\ \cline{2-6} 
                          & 5                    & $338$                   & $62.2720$                     & \multicolumn{1}{c|}{$6.59649 \times 10^{-3}$} & $3.91198 \times 10^{-4}$ \\ \cline{2-6} 
                          & 6                    & $411$                   & $99.5202$                     & \multicolumn{1}{c|}{$9.07259 \times 10^{-3}$} & $5.44030 \times 10^{-4}$ \\ \cline{2-6} 
                          & 7                    & $483$                   & $135.533$                     & \multicolumn{1}{c|}{$1.19532 \times 10^{-2}$} & $7.13111 \times 10^{-4}$ \\ \cline{2-6} 
                          & 8                    & $545$                   & $147.429$                     & \multicolumn{1}{c|}{$1.55013 \times 10^{-2}$} & $8.82928 \times 10^{-4}$ \\ \cline{2-6} 
                          & 9                    & $617$                   & $195.847$                     & \multicolumn{1}{c|}{$1.87248 \times 10^{-2}$} & $1.06232 \times 10^{-3}$ \\ \cline{2-6} 
                          & 10                   & $703$                   & $236.475$                     & \multicolumn{1}{c|}{$2.22630 \times 10^{-2}$} & $1.33978 \times 10^{-3}$ \\ \hline
    built-in                   & ——                   & ——                      & $37.1714$                     & \multicolumn{1}{c|}{——}                       & ——                       \\ \hline
    \end{tabular}%
    }
\end{table*}

\begin{table*}[]
    \centering
    \caption{Distributed Linear Least Squares (Complete Graph)}
    \label{tab:leastsquare_complete}
    \resizebox{0.9\textwidth}{!}{%
    \begin{tabular}{|c|c|c|c|cc|}
    \hline
    \multirow{2}{*}{Methods}   & \multirow{2}{*}{N} & \multirow{2}{*}{Average steps} & \multirow{2}{*}{Average running time/s} & \multicolumn{2}{c|}{Average error}                                                \\ \cline{5-6} 
                          &                      &                         &                               & \multicolumn{1}{c|}{$L_2$ norm}                  & $L_\infty$ norm             \\ \hline
    \multirow{9}{*}{PPCM} & 2                    & $30$                    & $3.00541$                     & \multicolumn{1}{c|}{$7.50367 \times 10^{-4}$} & $4.32202 \times 10^{-5}$ \\ \cline{2-6} 
                          & 3                    & $28$                    & $6.43045$                     & \multicolumn{1}{c|}{$7.13905 \times 10^{-4}$} & $4.45913 \times 10^{-5}$ \\ \cline{2-6} 
                          & 4                    & $34$                    & $9.21670$                     & \multicolumn{1}{c|}{$2.40951 \times 10^{-4}$} & $1.43743 \times 10^{-5}$ \\ \cline{2-6} 
                          & 5                    & $41$                    & $15.9335$                     & \multicolumn{1}{c|}{$1.18730 \times 10^{-4}$} & $7.47532 \times 10^{-6}$ \\ \cline{2-6} 
                          & 6                    & $52$                    & $26.6574$                     & \multicolumn{1}{c|}{$2.70836 \times 10^{-4}$} & $1.81067 \times 10^{-5}$ \\ \cline{2-6} 
                          & 7                    & $67$                    & $34.7499$                     & \multicolumn{1}{c|}{$2.74500 \times 10^{-4}$} & $1.60976 \times 10^{-5}$ \\ \cline{2-6} 
                          & 8                    & $69$                    & $36.9573$                     & \multicolumn{1}{c|}{$2.12387 \times 10^{-4}$} & $1.32010 \times 10^{-5}$ \\ \cline{2-6} 
                          & 9                    & $74$                    & $47.6610$                     & \multicolumn{1}{c|}{$2.34231 \times 10^{-4}$} & $1.45633 \times 10^{-5}$ \\ \cline{2-6} 
                          & 10                   & $107$                   & $71.9016$                     & \multicolumn{1}{c|}{$2.97032 \times 10^{-4}$} & $1.77270 \times 10^{-5}$ \\ \hline
    \multirow{9}{*}{WAGM} & 2                    & $216$                   & $13.2036$                     & \multicolumn{1}{c|}{$3.61078 \times 10^{-3}$} & $2.13656 \times 10^{-4}$ \\ \cline{2-6} 
                          & 3                    & $249$                   & $20.8617$                     & \multicolumn{1}{c|}{$2.77573 \times 10^{-3}$} & $1.61357 \times 10^{-4}$ \\ \cline{2-6} 
                          & 4                    & $241$                   & $38.7346$                    & \multicolumn{1}{c|}{$2.71751 \times 10^{-3}$} & $1.70948 \times 10^{-4}$ \\ \cline{2-6} 
                          & 5                    & $227$                   & $51.3783$                     & \multicolumn{1}{c|}{$2.70917 \times 10^{-3}$} & $1.53580 \times 10^{-4}$ \\ \cline{2-6} 
                          & 6                    & $213$                   & $58.2196$                     & \multicolumn{1}{c|}{$2.99082 \times 10^{-3}$} & $1.76974 \times 10^{-4}$ \\ \cline{2-6} 
                          & 7                    & $197$                   & $56.8410$                     & \multicolumn{1}{c|}{$3.99842 \times 10^{-3}$} & $2.27264 \times 10^{-4}$ \\ \cline{2-6} 
                          & 8                    & $201$                   & $66.0159$                     & \multicolumn{1}{c|}{$5.52460 \times 10^{-3}$} & $3.18646 \times 10^{-4}$ \\ \cline{2-6} 
                          & 9                    & $218$                   & $80.1341$                     & \multicolumn{1}{c|}{$7.59362 \times 10^{-3}$} & $4.35835 \times 10^{-4}$ \\ \cline{2-6} 
                          & 10                   & $234$                   & $86.3198$                     & \multicolumn{1}{c|}{$1.01300 \times 10^{-2}$} & $5.78401\times 10^{-4}$  \\ \hline
    built-in                   & ——                   & ——                      & $37.1714$                     & \multicolumn{1}{c|}{——}                       & ——                       \\ \hline
    \end{tabular}%
    }
\end{table*}

\begin{table*}[]
    \centering
    \caption{Distributed Logistic Regression ($N=5$, Complete Graph)}
    \label{tab:logistic}
    \resizebox{0.95\textwidth}{!}{
    \begin{tabular}{|c|c|c|cc|}
        \hline
        \multirow{2}{*}{Methods} & \multirow{2}{*}{Average steps} & \multirow{2}{*}{Average running time/s} & \multicolumn{2}{c|}{Average error} \\
        \cline{4-5} 
        & & & \multicolumn{1}{c|}{$L_2$ norm} & $L_\infty$ norm \\
        \hline
        PPCM & $58$ & $0.231966$ & \multicolumn{1}{c|}{$3.39272 \times 10^{-4}$} & $1.57422 \times 10^{-4}$ \\
        \hline
        WAGM & $10990$ & $19.7866$ & \multicolumn{1}{c|}{$8.15459 \times 10^{-2}$} & $5.43103 \times 10^{-2}$ \\
        \hline
        built-in & —— & $6.30731$ & \multicolumn{1}{c|}{——} & —— \\ 
        \hline
    \end{tabular}
    }
\end{table*}

\begin{table*}[]
    \centering
    \caption{Distributed Support Vector Machine ($N=5$, Complete Graph)}
    \label{tab:svm}
    \resizebox{0.95\textwidth}{!}{
    \begin{tabular}{|c|c|c|cc|}
        \hline
        \multirow{2}{*}{Methods} & \multirow{2}{*}{Average steps} & \multirow{2}{*}{Average running time/s} & \multicolumn{2}{c|}{Average error} \\
        \cline{4-5} 
        & & & \multicolumn{1}{c|}{$L_2$ norm} & $L_\infty$ norm \\
        \hline
        PPCM & $43$ & $0.313544$ & \multicolumn{1}{c|}{$1.32018 \times 10^{-4}$} & $6.44231 \times 10^{-5}$ \\
        \hline
        WAGM & $163$ & $0.545714$ & \multicolumn{1}{c|}{$4.83286 \times 10^{-4}$} & $1.40330 \times 10^{-4}$ \\
        \hline
        built-in & —— & $2.09686$ & \multicolumn{1}{c|}{——} & —— \\ 
        \hline
    \end{tabular}
    }
\end{table*}

\section{Conclusion}
\par Building on the theoretical foundations of proximal point algorithms and projection contraction methods, this study innovatively proposes an Adaptive Projection-based Prediction-Correction Method (PPCM), specifically designed to address structured monotone variational inequality problems. This method leverages only the gradient information of the objective function for computation, significantly simplifying the computational process and enhancing the practical applicability of the approach. The selection of algorithm parameters is clear and concise, ensuring ease of operation and implementation while maintaining superior algorithm performance. The design of the adaptive adjustment criteria is both intuitive and convenient, and the theoretically established convergence properties provide a solid guarantee for the algorithm's stability and reliability. Moreover, careful enhancements to PPCM enable it to effectively tackle distributed consensus optimization problems, broadening its range of applications. The decentralized nature of PPCM is evidenced by its reliance on local information for network updates, further augmenting the algorithm's flexibility and autonomy. Through a series of numerical experiments, the exemplary efficiency and reliability of this method have been thoroughly demonstrated. Looking forward, we anticipate delving deeper into the challenges of distributed optimization and are committed to exploring advanced distributed optimization algorithms that support asynchronous iterations, with the aim of further advancing this field.

\bibliographystyle{ieeetr}
\bibliography{reference}

\end{document}